\documentclass[12pt,reqno]{amsart}
\usepackage{amsmath,amstext, amsthm,amsfonts,amssymb,amsthm,comment}
\usepackage{amssymb}
\usepackage{graphicx}
\usepackage{amscd}
\usepackage{color}
\usepackage{float}
\usepackage{graphics,amsmath,amssymb}
\usepackage{amsthm}

\usepackage{amsfonts}
\usepackage{latexsym}
\usepackage{epsf}
\usepackage{xifthen}
\usepackage{mathrsfs}
\usepackage{dsfont}
\usepackage{makecell}
\usepackage[FIGTOPCAP]{subfigure}
\usepackage[justification=centering,singlelinecheck=false]{caption}
\usepackage{amsmath}
\allowdisplaybreaks[4]
\usepackage{listings}
\usepackage{etoolbox}
\usepackage{fancyhdr}
\usepackage{pdflscape}
\usepackage[title,toc,titletoc]{appendix}
\usepackage{enumitem}
\usepackage[noadjust]{cite}
\usepackage{tikz}
\usepackage{hyperref} 

\newtheorem{theorem}{Theorem}[section]
\newtheorem{lemma}[theorem]{Lemma}
\newtheorem{corollary}[theorem]{Corollary}
\newtheorem{prop}[theorem]{Proposition}

\newtheorem{example}[]{Example}
\newtheorem{remark}[]{Remark}
\newtheorem{ques}[]{Question}

\renewcommand{\a}{\alpha}
\newcommand\qbin[3]{{\left[\begin{matrix} #1 \\ #2 \end{matrix} \right]_{#3}}}

\begin{document}

\begin{abstract}
Explicit positive multisums on base $q^5$ for the Rogers-Ramanujan 
sum function and its partial sum are given. 32 combinatorial interpretations for these 
identities are given. The classical mod 5 Rogers-Ramanujan identities are 
recovered as a limiting case. These results are equivalent to a 
quintic basic hypergeometric identity.
\end{abstract}

\title{Partial sums of Rogers-Ramanujan functions}
\author{Dennis Stanton}
\address{School of Mathematics, University of Minnesota, Minneapolis, MN 55455}
\email{stanton@math.umn.edu}
\date{November 18, 2026}
\maketitle

\section{Introduction} 
In 1999 Garrett, Ismail and Stanton \cite[Theorem 7.1]{GIS} gave 
a quintic transformation for a basic hypergeometric series
\begin{equation}
\label{quintic}
\sum_{n=0}^\infty \frac{q^{n^2} x^n}{(q;q)_n}=
\frac{(x^2q^5,xq^3,x^2q^2;q^5)_\infty}{(xq;q)_\infty}
\ _3\phi_2
\left(\begin{matrix}
x,&xq,&xq^3\\
&x^2q^5,&x^2q^2
\end{matrix}
\ \bigg|\ q^5,xq^3\right).
\end{equation}
One consequence of \eqref{quintic} was the Rogers-Ramanujan identity 
\cite[p. 104]{And} for $1,4 \mod 5$ by choosing $x=1$
$$
\sum_{n=0}^\infty \frac{q^{n^2}}{(q;q)_n}=
\frac{1}{(q^1,q^4;q^5)_\infty}.
$$

One may ask if \eqref{quintic} contains more combinatorial information
than just its value at $x=1$. The left side is manifestly positive as power series in $x$ and $q$. 
Can the right side of \eqref{quintic} be recast in a similar way, on base $q^5$, 
and reveal new combinatorial results? In this paper we answer this question affirmatively, see
Theorems~\ref{main}, \ref{main23}, and Corollaries~\ref{maincor} and \ref{maincor23}. 
The combinatorial interpretations for the partial sums in terms of $\mod 5$ partitions are 
given in Theorems~\ref{main2} and \ref{main2-23}.  

We use the standard notation found in \cite{And} and \cite{gasper}. We also let
$\ell(\lambda)$ denote the number of parts of a partition $\lambda$.

\section{The quintic transformation rewritten}

One way to rewrite the quintic transformation, upon dividing by $1-x$,  and 
cancelling products, is
\begin{prop}
\label{prop1}
\begin{equation}
\frac{1}{1-x}\sum_{j=0}^\infty \frac{q^{j^2}x^j}{(q;q)_j}=
\sum_{k=0}^\infty
\frac{(x^2q^{5+5k};q^5)_\infty}{ (xq^{1+5k};q^5)_\infty(xq^{4};q^5)_\infty}
\frac{(x^2q^{2+5k};q^5)_\infty}{ (xq^{5k};q^5)_\infty(xq^{2};q^5)_\infty}
\frac{(xq^3;q^5)_k}{(q^5;q^5)_k} (xq^3)^k
\end{equation}
\end{prop}

The right side is nearly positive. 
The first quotient has a combinatorial interpretation, positive, using Theorem~\ref{easy2} 
and parts which are 1,4 mod 5. 
The second quotient has a combinatorial interpretation, positive, again using Theorem~\ref{easy2} 
and parts which are 0,2 mod 5. 
The third quotient has a combinatorial interpretation, not all positive, 
the parts which are 3 mod 5. This due to the $(xq^3;q^5)_k$ 
term in the summation.

The rational function identity for the coefficient of $x^n$ 
equivalent to Proposition~\ref{prop1} is
\begin{equation}
\label{bigsum}
\begin{aligned}
\sum_{j=0}^n \frac{q^{j^2}}{(q;q)_j}=&
\sum_{\substack{a,b,c,d,k,s \ge 0\\a+b+c+d+k+s= n}}^n
\frac{q^{5ab} q^{4a}q^{b(5k+1)}}{(q^5;q^5)_a(q^5;q^5)_b} 
\frac{q^{5cd} q^{5kc}q^{2d}}{(q^5;q^5)_c(q^5;q^5)_{d}} \frac{q^{3k} }{(q^5;q^5)_k}
\\
\times &
\qbin{k}{s}{q^5} q^{5\binom{s}{2}}(-1)^s q^{3s}.
\end{aligned}
\end{equation}
Note that the only negative terms in \eqref{bigsum} arise from $(-1)^s$.

\begin{ques} Is there an involution for the right side, which is the generating function for
certain 5-tuples $(\lambda_1,\lambda_4,\lambda_0,\lambda_2,\lambda_3)$,
where $\lambda_i$ has parts which are $i\mod 5,$ 
$\lambda_3$ is marked, and whose fixed points are all positive?
\end{ques}

ChatGPT6.0 answered Question 1 via a weight-preserving involution, whose fixed point set 
has $s=0$, also inserts a $q^{5ad}$, and thus has all terms are positive.  The result is Theorem~\ref{main}.
We give an independent (human) proof of Theorem~\ref{main} in Lemma~\ref{mainlem} of \S4 
using generating functions, and apply 
it for the $2,3\mod 5$ case in \S6. Lemma~\ref{mainlem} appears in \cite[Corollary 3.4, Theorem 3.5]{ISV} as 
a generating function for crossings in inhomogeneous matchings in a 5-partite graph.

\section{Main results}

Here are the Main Results: Theorem~\ref{main} and Corollary~\ref{maincor}, 
for the mod 5 expansion of the Rogers-Ramanujan sum and its combinatorics, along with 
similar results for the partial sum 
$$
\sum_{n=0}^N \frac{q^{n^2}}{(q;q)_n}
$$ 
in Theorem~\ref{main2}. 

\begin{theorem}[Main Theorem]
\label{main}
Let $Q=q^5. $ Then 
$$
\frac{1}{1-x}\sum_{n=0}^\infty \frac{q^{n^2} x^n}{(q;q)_n}=
\sum_{n_0,n_1,n_2,n_3,n_4\ge 0} \frac{x^{n_0+n_1+n_2+n_3+n_4} q^{E(n)}}{\prod_{i=0}^4 (Q;Q)_{n_i}}.
$$
where
$$
E(n)=\sum_{i=0}^4 in_i +5\sum_{i=0}^4 n_in_{i+2}.
$$
The indices are taken modulo 5.
\end{theorem}

\begin{corollary} 
\label{maincor}
The right side of Theorem~\ref{main} is the generating function 
for all 5-tuples of partitions, $(\lambda_1,\lambda_4,\lambda_0,\lambda_2,\lambda_3)$ where 
\begin{enumerate}
\item the parts of each $\lambda_i$ are $i\mod 5$ and $\lambda_i$ has $n_i$ parts,
\item every part of $\lambda_0 \ge 5(\ell(\lambda_2)+\ell(\lambda_3)),$
\item every part of $\lambda_1 \ge 1+5(\ell(\lambda_4)+\ell(\lambda_3)),$
\item every part of $\lambda_2 \ge 2+5(\ell(\lambda_4)).$
\end{enumerate}
If $\lambda_i=\varnothing$ for $0\le i\le 2$, then the conditions (2)-(4) are vacuous.
\end{corollary}

Note: $0$ is an allowed part of $\lambda_0.$

There are 32 variations of Corollary~\ref{maincor}  
which contain 32 different combinatorial statements.
These depend upon the two choices for interpreting $q^{5n_i n_{i+2}}$, $0\le i\le 4$:
either $n_i$ parts each at least $5n_{i+2}$, or $n_{i+2}$ parts each at least $5n_i.$
The same may be said about Theorem~\ref{main2}.

By equating coefficients of $x^N$ in Theorem~\ref{main} we have the result for 
a mod 5 interpretation of the partial sums of the Rogers-Ramanujan identities.

\begin{theorem}[Second Main Theorem]
\label{main2}
Let $Q=q^5. $ Then 
$$
\sum_{n=0}^N \frac{q^{n^2}}{(q;q)_n}=
\sum_{\substack{n_0,n_1,n_2,n_3,n_4\ge 0\\n_0+n_1+n_2+n_3+n_4=N}} 
\frac{q^{E(n)}}{\prod_{i=0}^4 (Q;Q)_{n_i}}.
$$
The left side is the generating function for difference two partitions with at most $N$ parts. 
The right side is the generating function for 5-tuples of 
partitions given by Corollary~\ref{maincor}, with the additional condition that
the 5-tuple has a total of $N$ parts.
\end{theorem}

\begin{example} 
\label{Exam1}
Let $N=2$ in Theorem~\ref{main2}. The identity is
$$
\begin{aligned}
1+&\frac{q^1}{(1-q)}+\frac{q^4}{(1-q)(1-q^2)}\\
=&\frac{1+q^2+q^4+q^6+q^8}{(1-q^5)(1-q^{10})}
+\frac{q^1 + q^7 + q^8 + q^4 + q^3 + q^9 + q^{10} + q^5 + q^{11} + q^7}{(1-q^5)(1-q^5)}.
\end{aligned}
$$
The 15 5-tuples of partitions $(\lambda_1,\lambda_4,\lambda_0,\lambda_2,\lambda_3)$ 
with exactly two parts, are, in respective order,
$$
\begin{aligned}
&(\varnothing,\varnothing,ab,\varnothing,\varnothing), 
(ab,\varnothing,\varnothing,\varnothing,\varnothing),
(\varnothing,\varnothing,\varnothing,ab,\varnothing),
(\varnothing,\varnothing,\varnothing,\varnothing,ab),
(\varnothing,ab,\varnothing,\varnothing,\varnothing)\\
&(a,\varnothing,b,\varnothing,\varnothing),
(\varnothing,a,\varnothing,\varnothing,b), 
(\varnothing,\varnothing,a,\varnothing,b), a\ge 5,
(\varnothing,a ,b, \varnothing,\varnothing),
(a,\varnothing,\varnothing,b,\varnothing)\\
&(a,\varnothing,\varnothing,\varnothing,b),a\ge 6,
(a,b,\varnothing,\varnothing,\varnothing), a\ge 6,
(\varnothing, \varnothing,\varnothing,a,b),
(\varnothing,a,\varnothing,b,\varnothing), b\ge 7,\\
&(\varnothing,\varnothing,a,b,\varnothing), a\ge 5.
\end{aligned}
$$
\end{example}

\begin{example} Take $N=2$, and one of remaining 31 combinatorial interpretations,
\begin{enumerate}
\item the parts of each $\lambda_i$ are $i\mod 5$ and $\lambda_i$ has $n_i$ parts,
\item every part of $\lambda_0 \ge 5(\ell(\lambda_2)),$
\item every part of $\lambda_1 \ge 1+5(\ell(\lambda_3)),$
\item every part of $\lambda_2 \ge 2+5(\ell(\lambda_4)),$
\item every part of $\lambda_3 \ge 3+5(\ell(\lambda_0)).$
\item every part of $\lambda_4 \ge 4+5(\ell(\lambda_1)).$\end{enumerate}
The 15 5-tuples of partitions $(\lambda_1,\lambda_4,\lambda_0,\lambda_2,\lambda_3)$ 
with exactly two parts, are, in respective order,
$$
\begin{aligned}
&(\varnothing,\varnothing,ab,\varnothing,\varnothing), 
(ab,\varnothing,\varnothing,\varnothing,\varnothing),
(\varnothing,\varnothing,\varnothing,ab,\varnothing),
(\varnothing,\varnothing,\varnothing,\varnothing,ab),
(\varnothing,ab,\varnothing,\varnothing,\varnothing)\\
&(a,\varnothing,b,\varnothing,\varnothing),
(\varnothing,a,\varnothing,\varnothing,b), 
(\varnothing,\varnothing,a,\varnothing,b), b\ge 8,
(\varnothing,a ,b, \varnothing,\varnothing),
(a,\varnothing,\varnothing,b,\varnothing)\\
&(a,\varnothing,\varnothing,\varnothing,b),a\ge 6,
(a,b,\varnothing,\varnothing,\varnothing), b\ge 9,
(\varnothing, \varnothing,\varnothing,a,b),
(\varnothing,a,\varnothing,b,\varnothing), b\ge 7,\\
&(\varnothing,\varnothing,a,b,\varnothing), a\ge 5.
\end{aligned}
$$
\end{example}

\begin{remark} The allowable sets of 5-tuples are not monotonically increasing in 
$N$, even if 0 parts in $\lambda_0$ are deleted. In  Example~\ref{Exam1}, if $N=2$, 
$(\varnothing,\varnothing,5,2,\varnothing)$ is allowed. But it is not allowed for any $N\ge 3$ by 
attaching 0's to $\lambda_0$,  because condition (2) 
in Corollary~\ref{maincor} says the parts of $\lambda_0$ must be at least 5.
\end{remark}

\begin{ques} Ideally the allowable 5-tuples for $N$ would all have 
$\lambda_2=\lambda_3=\lambda_0=\varnothing,$ and some conditions on 
$\lambda_1$ and $\lambda_4$ which depend upon $N$. We do not know such a bijection.
Perhaps one of the 32 possibilities for inequalities is more receptive to a bijection.
\end{ques}

\section{A $q$-series proof}

We directly prove Theorem~\ref{main} using $q$-series. 
It follows from Lemma~\ref{mainlem} below and the quintic transformation.

\begin{lemma} 
\label{mainlem}
We have
$$
\begin{aligned}
R=&\sum_{n_0,n_1,n_2,n_3,n_4\ge 0} \biggl(\prod_{i=0}^4 x_i^{n_i}\biggr)
\frac{ Q^{n_0n_2+n_1n_3+n_2n_4+n_3n_0+n_4n_1}}{\prod_{i=0}^4 (Q;Q)_{n_i}}\\
=&
\frac{(x_3,x_1x_4,x_0x_2;Q)_\infty}{(x_0,x_1,x_2,x_3, x_4;Q)_\infty}
\times \ _3\phi_2
\left(\begin{matrix}
x_1,&x_2x_4/x_3,&x_0\\
&x_1x_4,&x_0x_2
\end{matrix}
\ \bigg|\ Q;x_3\right).
\end{aligned}
$$
\end{lemma}

\begin{proof}
The right side $R$ of Theorem~\ref{main} is
$$
R=\sum_{n_0,n_1,n_2,n_3,n_4\ge 0} \biggl(\prod_{i=0}^4 x_i^{n_i}\biggr)
\frac{ Q^{n_0n_2+n_1n_3+n_2n_4+n_3n_0+n_4n_1}}{\prod_{i=0}^4 (Q;Q)_{n_i}}.
$$

The $n_0$ sum may be evaluated to 
$$
1/(x_0Q^{n_2+n_3};Q)_\infty= \frac{(x_0;Q)_{n_2+n_3}}{(x_0;Q)_\infty}.
$$
$$
R=\sum_{n_1,n_2,n_3,n_4\ge 0} \biggl(\prod_{i=1}^4 x_i^{n_i}\biggr)
\frac{ Q^{n_1n_3+n_2n_4+n_4n_1}}{\prod_{i=1}^4 (Q;Q)_{n_i}} 
\frac{(x_0;Q)_{n_2+n_3}}{(x_0;Q)_\infty}..
$$

The $n_1$ sum may be evaluated to obtain
$$
1/(x_1Q^{n_3+n_4};Q)_\infty= \frac{(x_1;Q)_{n_3+n_4}}{(x_1;Q)_\infty},
$$
$$
R= \frac{1}{(x_0,x_1;Q)_\infty} \sum_{n_2,n_3,n_4\ge 0} \biggl(x_2^{n_2}x_3^{n_3}x_4^{n_4}\biggr)
\frac{ Q^{n_2n_4}}{ (Q;Q)_{n_2}(Q;Q)_{n_3}(Q;Q)_{n_4}}(x_0;Q)_{n_2+n_3} (x_1;Q)_{n_3+n_4}
$$

The $n_4$ sum may be evaluated by the $Q$-binomial theorem to obtain
$$
 \frac{(x_1x_4Q^{n_2+n_3};Q)_\infty}{(x_4Q^{n_2};Q)_\infty},
$$
$$
\begin{aligned}
R= &\frac{(x_1x_4;Q)_\infty}{(x_0,x_1,x_4;Q)_\infty} \sum_{n_2,n_3\ge 0} \biggl(x_2^{n_2}x_3^{n_3}\biggr)
\frac{1}{ (Q;Q)_{n_2}(Q;Q)_{n_3}}
\frac{(x_1;Q)_{n_3} (x_0;Q)_{n_2+n_3}(x_4;Q)_{n_2}}{(x_1x_4;Q)_{n_2+n_3}}\\
R=&\frac{(x_1x_4;Q)_\infty}{(x_0,x_1,x_4;Q)_\infty}\sum_{n_3\ge 0} \frac{(x_1;Q)_{n_3}(x_0;Q)_{n_3}x_3^{n_3}}
{(Q;Q)_{n_3}(x_1x_4;Q)_{n_3}}\ _2\phi_1(x_0Q^{n_3},x_4;x_1x_4Q^{n_3};Q;x_2)
\end{aligned}
$$

By Heine's transformation \cite[(III.1)] {gasper} we have
$$
\ _2\phi_1(x_0Q^{n_3},x_4;x_1x_4Q^{n_3};Q;x_2) =
\frac{(x_0Q^{n_3},x_2x_4;Q)_\infty}{(x_1x_4Q^{n_3},x_2;Q)_\infty}
\ _2\phi_1(x_1x_4/x_0,x_2;x_2x_4;Q;x_0Q^{n_3})
$$
so
$$
\begin{aligned}
R= &\frac{(x_2x_4;Q)_\infty}{(x_1,x_2,x_4;Q)_\infty} \sum_{s,n_3\ge 0} 
\frac{(x_1x_4/x_0,x_2;Q)_s}{ (Q;Q)_{s}(x_2x_4;Q)_{s}}x_0^s
\frac{(x_1;Q)_{n_3} }{(Q;Q)_{n_3}}x_3^{n_3}(Q^s)^{n_3}\\
=&\frac{(x_2x_4;Q)_\infty}{(x_1,x_2,x_4;Q)_\infty} \sum_{s\ge 0} 
\frac{(x_1x_4/x_0,x_2;Q)_s}{ (Q;Q)_{s}(x_2x_4;Q)_{s}}x_0^s
\frac{(x_1x_3Q^s;Q)_\infty}{(x_3Q^s;Q)_\infty}\\
=&\frac{(x_1x_3, x_2x_4;Q)_\infty}{(x_1,x_2,x_3,x_4;Q)_\infty}  
\ _3\phi_2
\left(\begin{matrix}
x_1x_4/x_0,&x_2,&x_3\\
&x_2x_4,&x_1x_3
\end{matrix}
\ \bigg|\ Q;x_0\right).
\end{aligned}
$$

However, \cite[III.10]{gasper} implies
$$
\begin{aligned}
\ _3\phi_2
\left(\begin{matrix}
x_1x_4/x_0,&x_2,&x_3\\
&x_2x_4,&x_1x_3
\end{matrix}
\ \bigg|\ Q;x_0\right)
=&
\frac{(x_3,x_1x_4,x_0x_2;Q)_\infty}{(x_1x_3,x_2x_4,x_0;Q)_\infty}\\
\times &\ _3\phi_2
\left(\begin{matrix}
x_1,&x_2x_4/x_3,&x_0\\
&x_1x_4,&x_0x_2
\end{matrix}
\ \bigg|\ Q;x_3\right).
\end{aligned}
$$
so 
$$
R=\frac{(x_3,x_1x_4,x_0x_2;Q)_\infty}{(x_0,x_1,x_2,x_3,x_4;Q)_\infty}\\
 \ _3\phi_2
\left(\begin{matrix}
x_1,&x_2x_4/x_3,&x_0\\
&x_1x_4,&x_0x_2
\end{matrix}
\ \bigg|\ Q;x_3\right).
$$
which completes the proof.
\end{proof}

\begin{proof}[Proof of Theorem~\ref{main}]
Choose $x_i=xq^i$ in Lemma~\ref{mainlem}, and use the quintic transformation \eqref{quintic}.
\end{proof}

\section{Rogers-Ramanujan: $N\to\infty$}

In this section we explore the combinatorial interpretations of 
Theorem~\ref{main2} when $N\to\infty.$ The classical interpretation is that the 
left side counts partitions whose parts differ by at least 2, while the 
right side counts partitions into parts congruent to 1 or 4 modulo 5. 

We begin with a near classical case.

\begin{prop} 
\label{prop1}
Fix a positive integer $M$. Then the number of partitions of $M$, into 
parts which differ by at least 2, equals the number of 5-tuples of partitions 
$(\lambda_1,\lambda_4,\lambda_0,\lambda_2,\lambda_3)$ of $M$
such that
\begin{enumerate}
\item the parts of each $\lambda_i$ are $i\mod 5$,
\item $\lambda_2=\varnothing$ and $\lambda_3=\varnothing,$
\item every part of $\lambda_1 \ge 1+5\ell(\lambda_4).$
\end{enumerate}
\end{prop}

\begin{proof}
Fix $M\ge 1$ and consider which terms in Theorem~\ref{main2} could contribute a
$q^M$. Since
$$
E(n)\ge n_1+n_2+n_3+n_4+n_0(n_2+n_3)\ge n_0+n_1+n_2+n_3+n_4=N
$$
if $n_2+n_3\ge 1.$ 
So if $N\ge M+1$ the exponent $E(n)$  includes  
only terms beyond $M$.  So the only terms that contribute to  Theorem~\ref{main2} have 
$n_2=n_3=0$. In this case the only remaining condition in Theorem~\ref{main2} is
that every part of $\lambda_1$ is at least $1+5\ell(\lambda_4).$

We note that dropping the zero parts on $\lambda_0$ allows any number of parts, not just $N$.
\end{proof}
 
\begin{example} 
 \label{RR1exam}
 Let $M=10$, there are 6 difference two partitions. The allowable 5-tuples are
  for Proposition~\ref{prop1}
  $$
 \begin{aligned}
 &(1111111111,\varnothing,\varnothing,\varnothing,\varnothing),
 \quad(61111,\varnothing,\varnothing,\varnothing,\varnothing),
 \quad(6,4,\varnothing,\varnothing,\varnothing), \\
 &(11111,\varnothing,5,\varnothing,\varnothing), 
\quad (\varnothing,\varnothing,55,\varnothing,\varnothing),
 \quad(\varnothing,\varnothing,10,\varnothing,\varnothing), 
 \end{aligned}
 $$
 \end{example}

The condition on $\lambda_1$ gives the generating function for $(\lambda_1, \lambda_4)$ 
(see Example~\ref{exam1})
$$
\frac{(q^5;q^5)_\infty}{(q^1,q^4;q^5)_\infty}
$$
while the $\lambda_0$ generating function is 
 $$
 \frac{1}{(q^5;q^5)_\infty}.
 $$
 So the $N\to\infty$ limit of Theorem~\ref{main2} is
\begin{equation}
\label{14mod5}
 \sum_{n=0}^\infty \frac{q^{n^2}}{(q;q)_n}=  \frac{(q^5;q^5)_\infty}{(q^1,q^4;q^5)_\infty} \frac{1}{(q^5;q^5)_\infty}
 =  \frac{1}{(q^1,q^4;q^5)_\infty}
\end{equation}

This is the classical result when the parts are 1,4 modulo.

 \begin{remark} The last expression is not a cancellation, but a bijection between
 $(\lambda_1,\lambda_4,\lambda_0)$, where every part of $\lambda_1 \ge 1+5(\ell(\lambda_4)),$ 
 and $\lambda_0$ unrestricted, 
 and  $(\lambda_1,\lambda_4,\varnothing)$, $(\lambda_1,\lambda_4)$ unrestricted.
 \end{remark}
 
Proposition~\ref{tripbij} is a bijection which also proves this equality.

\begin{prop}
\label{tripbij}
There is an explicit bijection between triples 
$$
\cup_{k=0}^{min(s,t)} \{(\lambda_1,\lambda_4,\lambda_0): 
\ell(\lambda_1)=s-k, {\text{every part of }} 
\lambda_1\ge 5(t-k)+1,\ \ell(\lambda_4)=t-k, \ \ell(\lambda_0)=k\}, 
$$
and triples 
$$
\{(\lambda_1,\lambda_4,\varnothing): 
\ell(\lambda_1)=s,\quad \ell(\lambda_4)=t\}.
$$ 
\end{prop}

\begin{proof}
This follows from the identity
$$
\sum_{k=0}^{min(s,t)}\frac{q^{s-k}q^{5(s-k)(t-k)}}{(Q;Q)_{s-k}} \frac{q^{4(t-k)}}{(Q;Q)_{t-k}}
\frac{q^{5k}}{(Q;Q)_k}= \frac{q^s}{(Q;Q)_s}\frac{q^{4t}}{(Q;Q)_t}.
$$
Using
$$
\frac{1}{(Q;Q)_{s-k}(Q;Q)_k}=\frac{1}{(Q;Q)_s}
\qbin{s}{k}{Q}
$$
which has a simple bijection, we need a bijection for 
\begin{equation}
\label{Durf}
\sum_{k=0}^{min(s,t)} Q^{(s-k)(t-k)}
\qbin{s}{k}{Q} \frac{1}{(Q;Q)_{t-k}}
= \frac{1}{(Q;Q)_t}.
\end{equation}

A bijection for \eqref{Durf} is provided by a Durfee rectangle. 
The right side is the generating function for all partitions $\mu$ whose parts are at most $t$. 
Let $k$ be the smallest integer such that the cell $(s-k,t-k)$ lies in the Ferrers diagram of $\mu$.
Put $k=min(s,t)$ if no such cell exists. 
Thus $\mu$ has a Durfee rectangle of size $(s-k)(t-k)$.  In the Ferrers diagram of $\mu$, to the right 
of the Durfee rectangle is a partition which lies inside an $(s-k)\times k$ rectangle. Below 
the Durfee rectangle is any partition whose largest part is at most $t-k.$ These, 
along with the Durfee rectangle, have respective generating functions of
$$
\qbin{s}{k}{Q}, \quad  \frac{1}{(Q;Q)_{t-k}}, \quad {\text{and}}\quad Q^{(s-k)(t-k)},
$$
which proves \eqref{Durf}.
 \end{proof}

\section{The second Rogers-Ramanujan identity: $2,3 \mod 5$}

We can use Lemma~\ref{mainlem} for the results which are $2,3\mod 5.$ 

\begin{theorem}
\label{main23}
Let $Q=q^5. $ Then 
$$
\frac{1}{1-x}\sum_{n=0}^\infty \frac{q^{n^2+n} x^n}{(q;q)_n}=
\sum_{n_0,n_6,n_2,n_3,n_4\ge 0} \frac{x^{n_0+n_6+n_2+n_3+n_4} q^{E_{23}(n)}}{\prod_{i=0}^4 (Q;Q)_{n_i}}.
$$
where
$$
E_{23}(n)=6n_6+2n_2+3n_3+4n_4 +5(n_0n_6+n_6n_2+n_2n_3+n_3n_4+n_4n_0).
$$
The indices are taken modulo 5, for example $n_6=n_1.$
\end{theorem}

\begin{proof} In Lemma~\ref{mainlem} replace $(n_0,n_1,n_2,n_3,n_4)$ by 
$(m_0,m_1,m_2,m_3,m_4).$ Then to obtain Theorem~\ref{main23}, put 
$$
m_0=n_0,\quad
m_2=n_6, \quad m_4=n_2, \quad m_1=n_3, \quad m_3=n_4,
$$
and
$$
x_0=x,\quad x_2=xq^6, \quad x_4=xq^2,\quad x_1=xq^3 \quad x_3=xq^4.
$$
so the right side is
$$
\frac{(xq^4,x^2q^5,x^2q^6;Q)_\infty}{(x,xq^2,xq^3,xq^4,xq^6;Q)_\infty}
 \ _3\phi_2
\left(\begin{matrix}
xq^3,&xq^4,&x\\
&x^2q^5,&x^2q^6
\end{matrix}
\ \bigg|\ Q;xq^4\right).
$$

Again \cite[III.10]{gasper} implies
$$
\begin{aligned}
\ _3\phi_2
\left(\begin{matrix}
xq^3,&xq^4,&x\\
&x^2q^5,&x^2q^6
\end{matrix}
\ \bigg|\ Q;xq^4\right)
=&
\frac{(xq^4,x^2q^7,x^2q^4;Q)_\infty}{(x^2q^5,x^2q^6,xq^4;Q)_\infty}\\
\times &\ _3\phi_2
\left(\begin{matrix}
xq,&xq^2,&xq^4\\
&x^2q^7,&x^2q^4
\end{matrix}
\ \bigg|\ Q;xq^4\right).
\end{aligned}
$$

The right side is 
$$
\frac{(xq^4,x^2q^7,x^2q^4;Q)_\infty}{(x,xq^2,xq^3,xq^4,xq^6;Q)_\infty}
 \ _3\phi_2
\left(\begin{matrix}
xq,&xq^2,&xq^4\\
&x^2q^7,&x^2q^4
\end{matrix}
\ \bigg|\ Q;xq^4\right).
$$
This is the right side of \eqref{quintic}, with $x$ replaced by $xq$, times $(1-xq)/(1-x).$
Thus we replace $x$ by $xq$  on the left side of \eqref{quintic}, and multiply by 
$(1-xq)/(1-x).$
\end{proof}

\begin{corollary} 
\label{maincor23}
The right side of Theorem~\ref{main} is the generating function 
for a 5-tuple of partitions, $(\lambda_6,\lambda_4,\lambda_0,\lambda_2,\lambda_3)$ where 
\begin{enumerate}
\item the parts of each $\lambda_i$ are $i\mod 5$ and $\lambda_i$ has $n_i$ parts, and $\lambda_6$ has no 1's,
\item every part of $\lambda_6 \ge 6+5(\ell(\lambda_2)),$
\item every part of $\lambda_0 \ge 5(\ell(\lambda_6)+\ell(\lambda_4)),$
\item every part of $\lambda_3 \ge 3+5(\ell(\lambda_4+\ell(\lambda_2)).$
\end{enumerate}
\end{corollary}

Note: $0$ is an allowed part of $\lambda_0.$

\begin{theorem}
\label{main2-23}
Let $Q=q^5. $ Then 
$$
\sum_{n=0}^N \frac{q^{n^2+n}}{(q;q)_n}=
\sum_{\substack{n_0,n_6,n_2,n_3,n_4\ge 0\\n_0+n_6+n_2+n_3+n_4=N}} 
\frac{q^{E_{23}(n)}}{\prod_{i=0}^4 (Q;Q)_{n_i}}.
$$
Moreover the right side is the generating function for 5-tuples of 
partitions given by Corollary~\ref{maincor23}, with the additional condition that
the 5-tuple has a total of $N$ parts..
\end{theorem}

Again there are 32 variations of Theorem~\ref{main23}, 
Theorem~\ref{main2-23}, and Corollary~\ref{maincor23}..

As in \S4, when $N\to\infty$ in Theorem~\ref{main2-23}, $\lambda_6=\lambda_4=\varnothing$, and 
what remains are $\lambda_0, \lambda_3,\lambda_2$ whose generating functions are
\begin{equation}
\label{23mod5}
 \sum_{n=0}^\infty \frac{q^{n^2+n}}{(q;q)_n}=  \frac{(q^5;q^5)_\infty}{(q^2,q^3;q^5)_\infty} \frac{1}{(q^5;q^5)_\infty}
 =  \frac{1}{(q^2,q^3;q^5)_\infty}.
\end{equation}

\section{Tiny results recalled}
We recall two ``tiny" results, which must have been known for many years. 
The author knew about them for 40 years, see \cite[(3.17)]{ISV}, but they may go back to Sylvester
or before. 

\begin{theorem} [Tiny theorem]
\label{easy1}
The quotient of infinite products $\frac{(xy;q)_\infty}{(x,y;q)_\infty}$ has non-negative 
power series coefficients in $q$, $x$, and $y$.
\end{theorem}

\begin{theorem} [Tiny theorem combinatorial]
\label{easy2}
In fact 
$$
\frac{(q^2xy;q)_\infty}{(qx,qy;q)_\infty}= \sum_{n,m=0}^\infty \frac{x^mq^{(n+1)m}}{(q)_m}
\frac{y^nq^{n}}{(q)_n},
$$ 
which is the generating function for pairs of partitions $(\lambda_1,\lambda_2),$ which
satisfy:
every part of $\lambda_1$ is at least $1+\ell(\lambda_2).$
\end{theorem}
\begin{proof} Let $\ell(\lambda_1)=m$ and $\ell(\lambda_2)=n$, 
so the smallest part of $\lambda_1$ is at least $n+1$.
\end{proof}

\begin{example} 
\label{exam1}
What does
$$
\frac{(x^2q^{5+5k};q^5)_\infty}{ (xq^{1+5k};q^5)_\infty(xq^{4};q^5)_\infty}
=\sum_{n,m=0}^\infty q^{5mn} \frac{(xq^{1+5k})^m}{(q^5;q^5)_m} 
 \frac{(xq^{4})^n}{(q^5;q^5)_n}
$$
count? 

The generating function for partitions $\lambda_4$ with parts $4\mod 5$, and 
exactly $n=\ell(\lambda_4)$ parts is
$$
\frac{x^n q^{4n}}{(q^5;q^5)_n}.
$$

The generating function for partitions $\lambda_1$ with parts $1\mod 5$, and 
exactly $m=\ell(\lambda_1)$ parts, each part at least
$5n+5k+1$  is
$$
q^{5mn} \frac{(xq^{1+5k})^m}{(q^5;q^5)_m}.
$$

It counts pairs of partitions $(\lambda_1, \lambda_4)$, where $\lambda_i$ has parts which are $i\mod 5$,
and each part of $\lambda_1$ is at least $5(\ell(\lambda_4))+5k+1.$
\end{example}

\section{Remarks}

One may ask for a bijective proof of Theorem~\ref{main2}, which would 
include a bijective proof of the first Rogers-Ramanujan identity when $N\to\infty.$
The difference two side has a similar expansion by considering the columns (modulo 5) 
of the partition formed by the double staircase removal. The requested result would be a bijective proof of  
\begin{equation}
\label{finalbig}
\sum_{\substack{n_0,n_1,n_2,n_3,n_4\ge 0\\n_0+n_1+n_2+n_3+n_4\le N}} 
q^{(n_0+n_1+n_2+n_3+n_4)^2}\frac{q^{n_1+2n_2+3n_3+4n_4}}{\prod_{i=0}^4 (Q;Q)_{n_i}} 
=\sum_{\substack{n_0,n_1,n_2,n_3,n_4\ge 0\\n_0+n_1+n_2+n_3+n_4 = N}} 
\frac{q^{E(n)}}{\prod_{i=0}^4 (Q;Q)_{n_i}}.
\end{equation}

Theorem~\ref{main2} writes the partial sum of the Rogers-Ramanujan function as 
the generating function for 5-tuples 
of partitions $(\lambda_1,\lambda_4,\lambda_0,\lambda_2,\lambda_3).$ If $N=\infty$, 
when $\lambda_2=\lambda_3=\varnothing$, Proposition~\ref{tripbij} rewrites the first three 
entries as $(\tilde\lambda_1,\tilde\lambda_4,\varnothing)$, where 
$\tilde\lambda_1$ and $\tilde\lambda_4$ are unconstrained. If we now consider the restrictions on 
$\tilde\lambda_1$ and $\tilde\lambda_4$ imposed by $\lambda_2$ and $\lambda_3$, 
we could recast Theorem~\ref{main2} as the generating function for 
$(\hat\lambda_1,\hat\lambda_4,\varnothing,\lambda_2,\lambda_3).$ For a fixed $N$, 
we would have a set of partitions into parts $1,4\mod 5,$ and $2,3\mod 5.$ As $N\to\infty$, 
the parts which are $2,3\mod 5$ disappear.  One may ask if Kadell's injection \cite{Kad}
could propel this idea into an entirely $1,4\mod5$ universe. 

\section*{Acknowledgments} The author thanks Kathleen O'Hara for many helpful discussions.

\end{document}